\documentclass[12pt]{article}
\usepackage{geometry,amsmath, amssymb, amsthm,graphicx,version,float,multirow,pgf,soul,cancel,wrapfig}
\usepackage{hyperref,tcolorbox,xcolor}
\usepackage{subcaption,titlesec}
\usepackage[mathscr]{euscript}
\usepackage[export]{adjustbox}
\restylefloat{figure}
\titleformat{\section}
{\normalfont\bfseries}
{\thesection.}{0.5em}{}
 \titleformat{\subsection}
{\normalfont\bfseries}
{\thesubsection.}{0.5em}{}
\begin{document}
\begin{center}	
{\vspace*{4cm}  \bf  Method of variation of parameters revisited}\\
	 { \vspace*{1cm} Swarup Poria \footnote{swarup\_p@yahoo.com} and Aman Dhiman\footnote{amandhiman11@yahoo.com} }
	\vspace*{0.5cm} 
	\\ {\it Department of Applied Mathematics,\\University of Calcutta,\\92 APC Road, Kolkata-700009, India.}\\
	\vspace*{0.5cm}
\end{center}
	\begin{abstract} 
		The  method of variation of parameter (VOP)  for solving linear ordinary differential equation is revisited in this article. Historically, Lagrange and Euler explained the method of variation of parameter in the context of perturbation method. In this article, we explain the construction of particular solutions of a linear ordinary differential equation  in the light of linearly independent functions in a more systematic way. 
In addition, we have shown that if the time variation of the `constants'  contribute substantially to the velocity then also the solution remains invariant. VOP method for system of $n$ linear ODE is discussed.  Duhamel’s principle has also been studied in
reference to a system of $n$ linear ODE for completeness of this review. Finally, applications of VOP method for constructing Green's function is reported.

	\end{abstract}


\section{Introduction}
The method of variation of parameter (VOP) is a technique for transforming solutions of a linear homogeneous ordinary differential equations into a particular integral of  the corresponding  inhomogeneous system.
The method of variation of parameter (VOP) and method of undetermined coefficients are two very useful methods for determining particular integrals of a linear  ordinary differential equations (ODE). However, 
the method of undetermined coefficients has two inherent weaknesses that limits its wider application to linear equations. Firstly, the method of undetermined coefficients  is only applicable to linear ODE with constant coefficients and secondly, the inhomogeneous part of the ODE  must be of some special type. On the other hand, the method of variation of parameters is superior due to  no such restriction.

\begin{wrapfigure}{r}{0.5\linewidth}
	\centering
	\includegraphics[width=5cm]{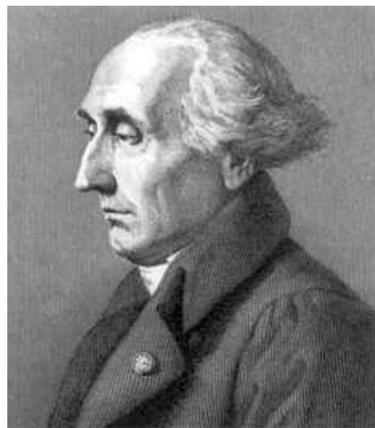}
	\caption*{Joseph Louis Lagrange}
\end{wrapfigure}
The method of variation of parameter was invented independently by Leonhard Euler (1748) and  by Joseph Louis Lagrange (1774). Although the method is famous for solving linear ODEs, it actually appeared in highly nonlinear context of celestial mechanics \cite{Newman}.  Euler and Lagrange were motivated to find solutions of the $n$-body problem of celestial mechanics with the help of the solution of two-body problem by converting the constants of motion or parameters into functions of time. This is the reason for the name of the method as  variation of parameters  or variation of constants.  It  is worth mentioning that Euler and Lagrange applied this method to nonlinear differential equations.  Lagrange gave the method of variation of parameters its final form during 1808-1810 .

In this review, a special emphasis is given on the way of constructing
particular integrals in the VOP method  in a  systematic way. The main contribution of  this review is to construct an elementary proof in support of the choice of particular integral in this method. In this technique there is a constraint on the time variation of the ``constants" which assumes that the time variability of the ``constants" does not contribute substantially to the velocity of the underlying dynamical equation represented by an ODE. In this article, we have generalized this constraint  by assuming that the time variation of the ``constants" can  contribute substantially to the velocity and show that the solution remains invariant under this generalization.  Duhamel's principle has also been discussed in context to a system of  $n$ linear ODE for completeness of this review. Finally, application of VOP method for constructing Green's function is reported.

\section{Construction of Particular Integral}
Consider a linear differential operator $L$ defined as follows,
$$L\equiv a_{n}(x)D^{n}+a_{(n-1)}(x)D^{n-1}+....+a_{1}(x)D+a_{0}(x)I, $$
$  ~\mbox{where}~D\equiv\frac{d}{dx},~I\equiv \mbox{identity operator}.$ The solutions of a linear homogeneous differential equation are called complementary functions. We define particular integral of a linear nonhomogeneous ODE as those functions which are not solutions of corresponding homogeneous linear ODE but satisfies the linear nonhomogeneous ODE. Let $y_c(x)$ and $y_p(x)$ be functions such that $Ly_c(x)=0  $ and $ Ly_p(x)=g(x), g(x) \ne 0. $  Then $y_c(x)$ is the complementary function and $y_p(x)$ is the particular integral. It is important to note that particular integral $y_p(x)$ and complementary functions $y_c(x)$ of any nonhomogeneous linear ODE are  linearly independent. The general solution of the ODE $ Ly(x)=g(x), g(x) \ne 0 $ is $ y=y_c(x)+ y_p(x).$\\

After application of same linear differential operator if one function produces zero value and other function produces a nonzero value then the two functions are linearly independent. It is a sufficient condition for linear independence of functions.  Notice that there can be $n$-linearly independent functions each of which produces a zero value when operated on a $n^{th}$ order homogeneous linear differential operator. This result was proved by Lagrange in 1765. Notice that this  beautiful property of linear differential operator  does not hold for nonlinear differential operators. In case of nonlinear differential operator  $L_1(f(x))=0$ but 
 $L_1(cf(x)) \ne 0$ for some nonzero scalar $c,$ although $f(x)$ and $cf(x)$ are linearly dependent functions.

Let us now consider  a set consisting  of two functions $\{ f_1(x), f_2(x) \} . $ If these two functions are linearly dependent on an interval then there exist constants $d_1$ and $d_2$ that are not both zero such that for every $x$ in the interval, $ d_1f_1(x)+ d_2f_2(x)=0.$ Therefore, without loss of generality  if we  assume that $d_1 \ne 0,$ then $ f_1(x)=-\frac{d_2}{d_1} f_2(x). $ 
Hence, if a set of two functions is linearly dependent, then one
function is  a constant multiple of the other. Conversely, if $f_1(x)=df_2(x)  $ for some constant $d$ then $ (-1).f_1(x) +d f_2(x)=0  $
 for every $x$ in the interval. Hence the set of functions is linearly dependent because at least one of the constants
$(-1)$ is not zero. Therefore, we can conclude that a set of two functions $f_1 (x)$ and $f_2 (x)$
is linearly independent when neither function is a constant multiple of the other on
the interval. Hence, if $y_1(x)$ and $y_2(x)$ are linearly independent then their quotient $\frac{y_2(x)}{y_1(x)}=u(x). $ As for example  the  functions $x$ and  $|x|$ is linearly independent on $(-\infty, \infty)$ as $\frac{|x|}{x}=\text{sgn(x)}.$

\subsection{First order linear ODE}
 Consider the first order linear ODE, 
\begin{eqnarray} \label{fir}
 \frac{dy}{dx} +p(x) y =q(x), ~ q(x) \ne 0    .
\end{eqnarray}
Let $y=d_1 y_c(x)$ be the general solution of the corresponding homogeneous equation 
$$ \frac{dy}{dx} +p(x) y =0    ,$$
where $d_1$ is an arbitrary constant.
Then the particular integral $y_p(x)$ is  a function such that $y_p(x)$ and $y_c(x)$ is  linearly independent.
Therefore,   $\frac{y_p(x)}{y_c(x)}=d_1(x)$ i.e., $y_p(x)=d_1(x) y_c(x),$ here $d_1(x)$ is an unknown function.
Substituting $y_p$ in \eqref{fir} we get
\begin{eqnarray}
 d'_1(x)y_c(x) & + & d_1(x) ( y_c'(x)+p(x) y_c(x))=q(x)  \nonumber \\
{\mbox{i.e.,}~~}~ d'_1(x) &=& \frac{q(x)}{y_c(x)} \nonumber 
\end{eqnarray}
Hence, $d_1(x) = \int \frac{q(x)}{y_c(x)} dx$ 
and the particular integral $y_p(x)=y_c(x)\int \frac{q(x)}{y_c(x)} dx .$ Therefore, the general solution is $ y(x)=y_c(x)+d_1(x)y_c(x).   $\\

\subsection{Second order linear ODE}

Consider the second order linear ODE 
\begin{eqnarray}\label{sec}
 \frac{d^2y}{dx^2} +p_1(x) \frac{dy}{dx}+ p_2(x) y =q(x), ~ q(x) \ne 0    .
\end{eqnarray}
In order to determine the particular integral $y_p(x)$ of  a second order linear ODE  using method of variation parameters we consider $y_p(x)= c_1(x) y_{1}(x) + c_2(x) y_{2}(x)$, where $y_{1}$ and $y_{2}$ are solutions to the corresponding homogeneous equation to \eqref{sec}.   In this section we present an intuitive proof in support of the above mentioned choice.\\

\begin{tcolorbox}

{\bf Theorem:}  The particular integral for a second order linear inhomogeneous ODE is given as $c_1(x) y_{1}(x) + c_2(x) y_{2}(x). $

\end{tcolorbox}

{\bf Proof :}Let $y= y_{1}(x)$ and $y= y_{2}(x)$   be two linearly independent   solutions of the  homogeneous equation corresponding to equation \eqref{sec},
$$  \frac{d^2y}{dx^2} +p_1(x) \frac{dy}{dx}+ p_2(x) y =0 .$$
Then the particular integral of the  linear ODE \eqref{sec} must be a function which is linearly independent to the set $\{y_{1}(x), y_{2}(x)\}.$ Clearly, there exist non constant functions $c_1(x)$ and $c_2(x)$ such that  $c_1(x) y_{1}(x)$ is linearly independent to $y_{1}(x)$ and $c_2(x) y_{2}(x)$  is linearly independent to $y_{2}(x).$  Note that  $c_1(x) y_{1}(x)$ is not necessarily  linearly independent to $y_{2}(x)$ and similarly $c_2(x) y_{2}(x)$ is not necessarily  linearly independent to $y_{1}(x)$.   We claim that 
$c_1(x) y_{1}(x) + c_2(x) y_{2}(x)$ is a function which is linearly independent to the set $\{y_{1}(x), y_{2}(x)    \}.$ We shall prove it by contradiction. Let there exist  scalars $d_1,d_2,d_3$ not all zero such that 
\begin{eqnarray}
d_1 y_{1}(x) + d_2 y_{2}(x) + d_3 \{c_1(x) y_{1}(x) + c_2(x) y_{2}(x)\} & = &0 ,  \nonumber\\
\mbox{i.e}~~~\{d_1 +d_3 c_1(x)\} y_{1}(x) + \{ d_2 +d_3 c_2(x))\} y_{2}(x)  & = & 0.
\end{eqnarray}

Then, $d_1 +d_3 c_1(x)=0 $ and $  d_2 +d_3 c_2(x))=0 $ implies that $c_1(x)$ and $c_2(x)$ are constant functions, which is a contradiction. Hence the only possibility is $d_1=d_2=d_3=0$ and the set $\{y_{1}(x), y_{2}(x), c_1(x) y_{1}(x) + c_2(x) y_{2}(x)    \}$ is a linearly independent set of functions.
Hence the proof is complete.

Therefore, it is clear that for a second order linear ODE the particular integral will be given by the function $c_1(x) y_{1}(x) + c_2(x) y_{2}(x).   $
One can  explain this form of particular integral  in the following way that the part of particular integral linearly independent to $y_{1}(x)$ is given by $c_1(x) y_{1}(x)$ and that of $y_{2}(x)$ is given by $c_2(x) y_{2}(x)$. The superposition of these two parts gives us the complete particular integral.  This result can be generalised for an $n^{th}$ order ODE in a straight forward manner.

\begin{tcolorbox}
	
{\bf Note :} Reduction of order (Euler 1753, Lagrange 1760)\cite{Anders} is a method which is applicable for finding general solution of any linear differential equation.
This reduction of order  method \cite{Reduce} converts  any linear differential equation to another linear differential equation of lower order provided at least one nontrivial solution of the ODE is known.  The general solution of  the original linear ODE can be obtained by using the general solutions of the lower-order equations.  Variation of parameter method can be viewed as a brilliant improvement of the reduction of order method for solving nonhomogeneous linear ODE.

\end{tcolorbox}

\section {Solution}

Let the particular intergral of (2) is given by
$$ y_p(x)= c_1(x) y_{1}(x) + c_2(x) y_{2}(x). $$

Let $y_1(x)$ and $y_2(x)$ be solutions to the homogeneous equation corresponding to \eqref{sec}. This gives 
\begin{align}\label{solu}
\begin{split}
y''_1+p_1(x)y'_1&+p_2(x)y_1=0,\\
y''_2+p_1(x)y'_2&+p_2(x)y_2=0.
\end{split}
\end{align}

By superposition principle $y(x)=c_1y_1(x)+c_2y_2(x)$ is also a solution for the  homogeneous equation corresponding to \eqref{sec}. The method of variation parameters assumes 
\begin{equation}\label{super}
y(x)=c_1(x)y_1(x)+c_2(x)y_2(x),
\end{equation}
to be the particular solution of equation \eqref{sec}, where $c_1(x),c_2(x)$ are to be determined. We find $y',y''$ of \eqref{super} to use in \eqref{sec}
\begin{align}\label{deriv}
\begin{split}
y'=c'_1(x)y_1(x)&+c'_2(x)y_2(x)+c_1(x)y'_1(x)+c_2(x)y'_2(x),\\
y''=c''_1(x)y_1(x)&+2c'_1(x)y'_1(x)+2c'_2(x)y'_2(x)+c''_2(x)y_2(x)+c_1(x)y''_1(x)+c_2(x)y''_2(x).
\end{split}
\end{align}
After substituting $y,y'$ and $y''$ in \eqref{sec} and doing a bit  of rearrangement we obtain,
\begin{align}\label{fin}
\begin{split}
\frac{d}{dx}\Big\{c'_1(x)y_1(x)+c'_2(x)y_2(x)\Big\}+p_1(x)\Big\{c'_1(x)&y_1(x)+c'_2(x)y_2(x)\Big\}\\
+&y'_1(x)c'_1(x)+y'_2(x)c'_2(x)=q(x).
\end{split}
\end{align}
In method of variation parameters for easier computation we choose 
\begin{equation}\label{1}
c'_1(x)y_1(x)+c'_2(x)y_2(x)=0,
\end{equation}
then the equation \eqref{fin} reduces to a first order ODE with two unknowns $c_1'(x)$ and $c_2'(x)$ of the following form,
\begin{equation}\label{2}
y'_1(x)c'_1(x)+y'_2(x)c'_2(x)=q(x).
\end{equation}
In the well known variation of parameter method then the work is to determine the parameters $c_1(x)$ and $c_2(x)$ by solving \eqref{1} and \eqref{2}. 

 Our motivation is to show that instead of setting  $c'_1(x)y_1(x)+c'_2(x)y_2(x)=0$ if we choose it as an arbitrary constant or more generally as a arbitrarily selected differentiable function  $A(x)$ then also one can  determine the unknown functions $c_1(x)$ and $c_2(x)$  to get the same solution $y(x)$.   In general we choose,

\begin{equation}\label{new1}
c'_1(x)y_1(x)+c'_2(x)y_2(x)=A(x).
\end{equation}
Substituting this to equation \eqref{fin}, we obtain,
\begin{equation}\label{new2}
y'_1(x)c'_1(x)+y'_2(x)c'_2(x)=q(x)-\left(A'(x)+p_1(x)A\right).
\end{equation} 
 On solving we get, 
$$c'_1(x)=\frac{\Big(q(x)-A'(x)-A(x)\Big)y_2(x)-A(x)y'_2(x)}{W(y_1(x),y_2(x))},$$ $$c'_2(x)=\frac{\Big(q(x)-A'(x)-p_1(x)A(x)\Big)y_1(x)-A(x)y'_1(x)}{-W(y_1(x),y_2(x))},$$
on integrating
$$c_1(x)=\int\frac{(q(x)-A'(x)-A(x))y_2(x)-A(x)y'_2(x)}{W(y_1(x),y_2(x))}dx,$$
$$c_2(x)=\int \frac{\Big(q(x)-A'(x)-p_1(x)A(x)\Big)y_1(x)-A(x)y'_1(x)}{-W(y_1(x),y_2(x))}dx.$$
Here $W(y_1(x),y_2(x))=y_1(x)y'_2(x)-y'_1(x)y_2(x)$ is the Wronskian for the solution. The particular solution is given as 
\begin{align} \label{finalbv}
y_p(x)&=\int\frac{y_1(x)y_2(s)-y_2(x)y_1(s)}{W(y_1(s),y_2(s))}q(s)ds-\int \frac{y_1(x)y_2(s)-y_2(x)y_1(s)}{W(y_1(s),y_2(s))}A'(s)ds\nonumber\\
-\int &\frac{y_1(x)y_2(s)-y_2(x)y_1(s)}{W(y_1(s),y_2(s))}p(s)A(s)ds-\int \frac{y_1(x)y'_2(s)-y_2(x)y'_1(s)}{W(y_1(s),y_2(s))}A(s)ds.
\end{align} 

\begin{tcolorbox}
{\bf Note :} Clearly the $q(x)$ independent integrals of \eqref{finalbv} will not be particular integral. Hence the total contribution from $A(x)$ and $A'(x)  $ dependent integral terms must be zero to the particular integral. It is an open problem to show that the particular solution $y_p(x)$ is completely independent of the choice of $A(x)$ in a more prominent way. 
\end{tcolorbox}

Therefore, the complete solution of \eqref{sec} is given as
\begin{align*}
y(x) =& y_{c}(x)+y_{p}(x)\\ 
     =&c_1 y_{1}(x) + c_2y_{2}(x)+\int\frac{y_1(x)y_2(s)-y_2(x)y_1(s)}{W(y_1(s),y_2(s))}q(s)ds.
\end{align*}

\section*{Example 1}
Consider the following differential equation
\begin{equation}\label{eg1}
y''-y'-2y=2e^{-x}.
\end{equation}
The auxiliary equation for \eqref{eg1} is $m^2-m-2=0$ which on solving gives to independent functions $e^{2x}$ and $e^{-x}$ as solutions to homogeneous equation corresponding to \eqref{eg1}. Here the complementary function is given as $y_c(x)=c_1 e^{2x} +c_2 e^{-x}$ and applying method of variation parameters, the particular solution is assumed as $$y_p(x)=c_1(x) e^{2x} +c_2(x) e^{-x}.$$ Now making use of \eqref{new1} and \eqref{new2}, we get
\begin{align}
c'_1(x) e^{2x} +c'_2(x) e^{-x}&=A(x),\label{egeq1}\\
2c'_1(x) e^{2x} -c'_2(x) e^{-x}&=2e^{t}+A(x)-A'(x).\label{egeq2}
\end{align}
Solving \eqref{egeq1} and \eqref{egeq2} for $c'_1(x)$ and $c'_2(x)$ we get
\begin{eqnarray}
c'_1(x)&=&\frac{2}{3}e^{-3x}+\frac{2}{3}A(x)e^{-2x}-\frac{1}{3}A'(x)e^{-2x},\label{c1}\\
c'_2(x)&=&-\frac{2}{3}+\frac{1}{3}A(x)e^x+\frac{1}{3}A'(x)e^x.\label{c2}
\end{eqnarray}
On integrating \eqref{c1} and \eqref{c2} we get
\begin{align*}
c_1=&\frac{2}{3}\int e^{-3x}+\frac{2}{3}\int A(x)e^{-2x}-\frac{1}{3}\int A'(x)e^{-2t}\\
   =&-\frac{2}{9} e^{-3x}+\frac{2}{3}\left[\frac{A(x)e^{-2x}}{-2}-\int \frac{A'(x)e^{-2x}}{-2}\right]-\frac{1}{3}\int A'(x)e^{-2x}\\
   =&-\frac{2}{9} e^{-3x}-\frac{A(x)e^{-2x}}{3}+\cancel{\frac{1}{3}\int A'(x)e^{-2x}}-\cancel {\frac{1}{3} \int A'(x)e^{-2x}}\\
   =&-\frac{2}{9} e^{-3x}-\frac{A(x)e^{-2x}}{3}.
\end{align*}
Similarly, on integrating \eqref{c2} we get 
\begin{equation}
c_2(x)=-\frac{2x}{3}+\frac{A(x)}{3}e^x.
\end{equation}
Using this the particular solution takes the form 
\begin{align}
y_p(x)=&\left(-\frac{2}{9}e^{-3x}-\frac{A(x)}{3}e^{-2x}\right)e^{2x}+\left(-\frac{2}{3}x+\frac{A(x)}{3}e^x\right)e^{-x}\\
       =&-\frac{2}{9}e^{-x}-\frac{2}{3}xe^x.
\end{align}
This is clearly independent of the choice of $A(x)$. 
\section*{Example 2}

Consider the following differential equation
\begin{equation}\label{ex2}
xy''-(x+1)y'+y=x^2,
\end{equation}
where it is given that $y_1(x)=e^x$ and $y_2(x)=1+x$ form the fundamental set of solutions for the homogeneous equation corresponding to \eqref{ex2}. Hence the complementary function is  $y_c(x)=c_1e^x+c_2(1+x).$ Applying the variation of parameters method, the particular solution is assumed as $y_p(x)= e^xc_1(x)+(1+x)c_2(x).$ Making use of \eqref{new1} and \eqref{new2}, with a special choice $A(x)=x^2$, we get
\begin{eqnarray}
e^xc'_1(x)+(x+1)c'_2(x)&=x^2,\label{egeqq1}\\
e^xc'_1(x)+c'_2(x)&=x^2.\label{egeqq2}
\end{eqnarray}
Solving \eqref{egeqq1} and \eqref{egeqq2} for $c'_1(x),c'_2(x)$ we get 
$$c'_1(x)=x^2e^{-x}~~,~~c'_2(x)=0$$  
on integrating, 
$$c_1(x)=-(x^2+2x+2)e^{-x}~~,~~c_2(x)=k$$
 Since the particular integral $y_p(x)$ should be free from arbitrary constants therefore we have to set  $k=0$ here. 

Therefore, the particular solution is given as 
$$y_p(x)=-(x^2+2x+2).$$
Finally, the complete solution is given as
\begin{equation*}
	y(x)=c_1e^x+c_2(1+x)-(x^2+2x+2).
\end{equation*}
An advantage of a particular choice of $A(x)$ is clear from the example 2.
\section{System of  ODEs}
Consider the system of linear differential equations in the following form,
\begin{equation}\label{sys}
x'(t)-P(t)x(t)=0,
\end{equation}
where $P_{n\times n}(t)$ is continuous on an interval I and $x_{n\times 1}$ is a column vector $[x_1(t),x_2(t),....x_n(t)]^t$. If $\phi_1(t),\phi_2(t),\phi_3(t),...,\phi_n(t)$ be $n$ linearly independent solutions of \eqref{sys} then the general solution $x(t)$ for \eqref{sys} is given by
$$x(t)=\sum_{i=1}^{n}c_{i}\phi_{i}(t)=\Phi(t)C,$$
where $C=(c_1,c_2,c_3,...,c_n)^t$ is a column vector and 

$$\Phi(t)=[\phi_1(t),\phi_2(t),\phi_3(t),...,\phi_n(t)]=\begin{pmatrix} 
\phi_{11} & \phi_{12} & \phi_{13} & .....& \phi_{1n} \\
\phi_{21} & \phi_{22} & \phi_{23} & .....& \phi_{2n} \\
\phi_{31} & \phi_{32} & \phi_{33} & .....& \phi_{3n} \\
\vdots & \vdots &\ddots  &\vdots  &\vdots\\
\phi_{n1} & \phi_{n2} & \phi_{n3} & .....& \phi_{nn} 
\end{pmatrix},$$
 is the fundamental matrix for \eqref{sys} on the interval I. In general the matrix $\Phi(t)$ is called a fundamental matrix   if its columns form a set of $n$-linearly independent solutions of \eqref{sys}. The column vector $C$ depends on the initial conditions $x(t_0)=x_0$. 
Substituting this initial condition  in the solution we obtain  $x_0= \Phi(t)C$ which gives $C= \Phi^{-1}(t_0)x_0 $ and the solution of the homogeneous equation can be written as  $x(t)=\Phi(t)\Phi^{-1}(t_0)x_0$. In the inhomogeneous case,
\begin{equation}\label{insys}
x'(t)-P(t)x(t)=b(t),
\end{equation}
where $b(t)$ is a column vector. 
According to the method of variation of parameter the particular  solution of \eqref{insys}  is of the form $x(t)=\Phi(t)C(t).$  Substituting  this in \eqref{insys}, we obtain
$$\Phi'(t) C(t)+\Phi(t) C'(t)=P(t)\Phi(t) C(t)+b(t).$$
 Since $\Phi'=P\Phi$, this gives $C'=\Phi^{-1}b$. Hence, $$C(t)=\int_{t_0}^{t}\Phi^{-1}(s)b(s)ds.$$
Therefore, the general solution   of \eqref{insys} is given by the superposition of solution of homogeneous part and the particular intergral as the following,
\begin{equation}\label{solusys}
 x(t)=\Phi(t)\Phi^{-1}(t_0)x_0+\int_{t_0}^{t}\Phi(t) \Phi^{-1}(s)b(s)ds.
\end{equation}

Obviously, if we solve the ODE \eqref{insys} with initial condition $x_0=0$ then we will get directly the  particular integral which is clear from \eqref{solusys}. This can be represented in the following manner,\\

\setlength\fboxrule{1pt}
{\centering
\begin{tabular}{ccc}\label{gen}
	\fcolorbox{black!50!black}{white}{$
	\begin{aligned}
	\textrm{Solution of}\quad \eqref{insys}\quad\textrm{with} \\ 
	\qquad \textrm{I.C.} \quad x(t_0)=x_0 
	\end{aligned}
	$} &
=
	\fcolorbox{black!50!black}{white}{$
		\begin{aligned}
	\textrm{Solution of} \quad\eqref{sys} \quad\textrm{with} \\ 
	\textrm{I.C.}  \quad x(t_0)=x_0 
		\end{aligned}
		$} &
	+
	\fcolorbox{black!50!black}{white}{$
		\begin{aligned}
	 \mbox{Solution of}\quad \eqref{insys} \mbox{ with} \\
	 \mbox{homog. I.C.} \quad x(t_0)=0 
		\end{aligned}
		$} \\
\end{tabular} 
}

Let $S_\tau^t=\Phi(t) \Phi^{-1}(\tau)$ be the solution operator for \eqref{sys}, then the solution of the IVP \eqref{insys} with $x(t_0)=x_0$ is given as 
$x(t)=S_{t_0}^tx_0+\int_{t_0}^{t}S_{s}^t(b(s))ds$. This is known as  Duhamel's principle \cite{washington}. This gives us a way to solve nonhomogeneous  linear differential equations, by superposition of solutions of corresponding homogeneous equation. In case of linear inhomogeneous ODEs  Duhamel's principle reduces to the method of variation of parameters. In general, Duhamel's principle is a  method for obtaining solutions to inhomogeneous linear evolution equations like the heat equation, wave equation etc. Duhamel principle allows to reduce the Cauchy problem for linear inhomogeneous partial differential equations to the Cauchy
problem for corresponding homogeneous equations.

\section{Green's function  }
George Green first published work on Green's function in 1828. Green's function has become a powerful tool for solving partial differential equations since then. In this article, we shall restrict ourself   to Green's functions for ordinary differential
equations. We will identify the Green's function for both initial
value and boundary value problems of nonhomogeneous second order linear differential equations \cite{Green} of the form 
\begin{equation} \label{diffe}
y''+p_1(x)y'+p_2(x)y=q(x), 
\end{equation}
on the  interval $x\in [a,b]$ using the method of variation of parameters.  It is important to note  that Green's function depends only on the solution functions  of homogeneous ODE (i.e. complementary function) but does not depend on the inhomogeneous  term $q(x)$.
\subsection{ Boundary Value Problem}
In case of boundary value problem our aim  is to find solution of \eqref{diffe} subject to boundary conditions $y(a)=y(b)=0$. Let $y_1(x)$ and $y_2(x)$ be two linearly independent solutions to the homogeneous equation corresponding to \eqref{diffe}. 
We want to express the solution of \eqref{diffe} in the following form,
\begin{equation}\label{green}
y(x)=\int_{a}^{b}G(x,s)q(s)ds,
\end{equation} 
where $G(x,s)$ is called the Green's function.  The  solution \eqref{finalbv} can be rearranged in the following form,

\begin{align*}
y(x)=&\int_{a}^{b}\frac{y_1(x)y_2(s)-y_2(x)y_1(s)}{W(y_1(s),y_2(s))}q(s)ds-\int_{a}^{b} \frac{y_1(x)y_2(s)-y_2(x)y_1(s)}{W(y_1(s),y_2(s))}A'(s)ds\\
&-\int_{a}^{b} \frac{y_1(x)y_2(s)-y_2(x)y_1(s)}{W(y_1(s),y_2(s))}p(s)A(s)ds-\int_{a}^{b} \frac{y_1(x)y'_2(s)-y_2(x)y'_1(s)}{W(y_1(s),y_2(s))}A(s)ds.
\end{align*}
Since the $q(x)$ independent integrals of \eqref{finalbv} will not be particular solution of the ODE, but it may produce terms linearly dependent to the complementary function of the ODE. Hence the total contribution from the last three  integral terms to the particular integral must be zero, which is consistent with the 
 particular case  $A(x)=0.$ In this case it reduces to an integral of the form  \eqref{green},where $$G(x,s)=\frac{y_1(x)y_2(s)-y_2(x)y_1(s)}{W(y_1(s),y_2(s))}.$$

\subsection{Initial Value Problem}
The problem is to find solution of equation \ref{gen} subject to $y(a)=y_0,\quad y'(a)=y'_0$. We make use Duhamel's principle to solve this problem by splitting \eqref{diffe} in two sub-problems as follows, 
\begin{equation}\label{ivp_sys}
y''+p_1(x)y'+p_2(x)y=0, \quad y_c(a)=y_0,\quad y'_c(a)=y'_0,\quad a<x<b,
\end{equation}
to obtain the complementary function and 
\begin{equation}\label{icp}
y''+p_1(x)y'+p_2(x)y=q(x), \quad y_p(a)=0,\quad y'_p(a)=0,\quad a<x<b,
\end{equation}   
to obtain  the particular solution.
The complementary function will be 
$$y_c(x)=c_1y_1(x)+c_2y_2(x),$$ 
where $c_1$ and $c_2$ are determined using
$$\begin{pmatrix} 
	y_1(a) & y_2(a) \\
	y'_1(a) & y'_2(a) 
\end{pmatrix}
\begin{pmatrix} 
c_1  \\
c_2  
\end{pmatrix}=
\begin{pmatrix} 
y_0 \\
y'_0 
\end{pmatrix}.  $$
For the particular solution we propose $y_p(x)$ as
$$y_p(x)=c_1(x)y_1(x)+c_2(x)y_2(x),$$
with the condition $c'_1(x)y_1(x)+c'_2(x)y_2(x)=A(x)$. Then, substituting in ODE \eqref{ivp_sys}, we can determine $c_1(x)$ and $c_2(x)$ by solving the system
$$\begin{pmatrix} 
y_1(a) & y_2(a) \\
y'_1(a) & y'_2(a) 
\end{pmatrix}
\begin{pmatrix} 
c'_1(x)  \\
c'_2(x)  
\end{pmatrix}=
\begin{pmatrix} 
0 \\
q(x) 
\end{pmatrix}.  $$
On solving for  $c'_1(x)$ and $c'_2(x)$ we get
\begin{align}\label{integ}
c'_1(x)=&-\frac{y_2(x)q(x)}{W(y_1(x),y_2(x))}+\frac{y'_2(x)A(x)}{W(y_1(x),y_2(x))},\\
c'_2(x)=&\frac{y_1(x)q(x)}{W(y_1(x),y_2(x))}-\frac{y'_1(x)A(x)}{W(y_1(x),y_2(x))}.
\end{align}  
to obtain $c_1(x)$ and $c_2(x)$ we need to integrate, for this we need the limits of integration. We know the initial conditions for $y_p$ by \eqref{icp}, these can be put as
$$\begin{pmatrix} 
y_1(a) & y_2(a) \\
y'_1(a) & y'_2(a) 
\end{pmatrix}
\begin{pmatrix} 
c_1(a)  \\
c_2(a)  
\end{pmatrix}=
\begin{pmatrix} 
0 \\
0
\end{pmatrix}.  $$
This on solving yields $c_1(a)=c_2(a)=0$. Furthermore, on performing the integration on \eqref{integ} we get
\begin{align*}
\int_{a}^{x} c'_1(s)ds=& c_1(x)= -\int_{a}^{x}\frac{y_2(s)}{W(y_1(s),y_2(s))}q(s)ds+\int_{a}^{x}\frac{y'_2(x)A(x)}{W(y_1(x),y_2(x))}ds,\\
\int_{a}^{x} c'_2(s)ds=& c_2(x)= \int_{a}^{x}\frac{y_1(s)}{W(y_1(s),y_2(s))}q(s)ds-\int_{a}^{x}\frac{y'_1(x)A(x)}{W(y_1(x),y_2(x))}ds,
\end{align*}
and then the particular solution
\begin{align}
&~~~~~~~~~~~~~~~~~~~~~~~~~y_p(x)=c_1(x)y_1(x)+c_2(x)y_2(x),\\
y_p(x)&=\int_{a}^{x}\frac{y_1(s)y_2(x)-y_1(x)y_2(s)}{W(y_1(s),W(y_2(s)))}q(s)ds+\int_{a}^{x}\frac{y_1(x)y'_2(s)-y'_1(s)y_2(x)}{W(y_1(s),W(y_2(s)))}A(s)ds.
\end{align}
Now with $A(x)=0$, rewriting the solution as 
$$y_p(x)=\int_{a}^{b}G(x,s)q(s)ds,$$
where
$$G(x,s)=\begin{cases} \frac{y_1(s)y_2(x)-y_1(x)y_2(s)}{W(y_1(s),W(y_2(s)))}  &\mbox{if } a<s<x \\ 
\qquad 0 & \mbox{if } x<s<b. \end{cases}$$

Furthermore, if the given boundary conditions are linearly independent, then the problem is well-defined and solvable using the method of variation parameters. Green's function can be determined for different cases of unmixed boundary conditions e.g. Dirichlet condition, Neumann condition and Robin conditions also for the mixed types like in periodic and anti-periodic conditions in a similar manner.

\begin{tcolorbox}
	
{\bf Note :} In case of equation  \eqref{insys} the solution with initial condition is given as follows
\begin{equation}
x(t)=\Phi(t)\Phi^{-1}(t_0)x_0+\int_{t_0}^{t} \underbrace{ \Phi(t) \Phi^{-1}(s)}_{G(t,s)}b(s)ds,
\end{equation}
where, 
$G(t,s)=\Phi(t) \Phi^{-1}(s)$ is the green's function here.

\end{tcolorbox}

\section{Conclusion  }

The VOP method is a very powerful technique for solving linear ODEs with  the help of complementary functions. Therefore, the method fails to solve inhomogeneous ODE  whenever the complementary functions cannot be determined. In this review, a novel explanation of the method for  construction of particular solutions  is given in the light of linearly independent functions in a more systematic way in contrast to the usual explanation via perturbation method. In the conventional VOP technique there is a constraint on the time variation of the ``constants"
which assumes that the time variability of the ``constants" does not contribute substantially to the velocity of the underlying dynamical equation represented by an ODE. In
this article, we have generalized this constraint by assuming that the time variation of the
“constants” can contribute substantially to the velocity and show that the general solution remains invariant under this generalization. Duhamel’s principle has also been discussed in
context to a system of $n$ linear ODE for completeness of this review. Construction of Green's function  through VOP  method is discussed.

We have discussed the method of variation of parameter for solving linear ODE only. However, general form of  this method can be used to solve nonlinear ODE \cite {Alekseev}, integro-differential equations \cite{Brunner} and even for solving  nonlinear functional differential equations \cite{Deo_and_Torres}. One can go through the contributions of Lakshmikantham \cite{Lakshmi} for understanding comprehensive applications of this method to differential equations. Deeper understanding of variation of parameters method is closely related to important topics like differential geometry, Lie symmetries, and the notions of reduction of order.



\end{document}